\documentclass{article}
\usepackage{amsmath,amssymb,amsthm,latexsym,graphicx}

\theoremstyle{plain}
\newtheorem{lemma}{Lemma}

\newtheorem{theorem}[lemma]{Theorem}
\newtheorem{corollary}[lemma]{Corollary}
\theoremstyle{definition}

\newtheorem{example}{Example}
\theoremstyle{remark}
\newtheorem{remark}{Remark}

\newcommand{\pdeg}{\mathop{\mathrm{pdeg}}}

\title{The Parametric Degree of a Rational Surface}

\author{Josef Schicho\thanks{
	The author was supported by the FWF in the frame of
		the reseach projects SFB 1303 and P15551} \\
	Johann Radon Institute \\
	Austrian Academy of Sciences}

\begin{document}

\maketitle

\begin{abstract}
The parametric degree of a rational surface is the degree of the polynomials
in the smallest possible proper parametrization. 
An example shows that the parametric degree 
is not a geometric but an arithmetic concept, in the sense that it
depends on the choice of the ground field. 
In this paper, we introduce two geometrical invariants of a rational surface,
namely level and keel. These two numbers govern the parametric degree
in the sense that there exist linear upper and lower bounds.
\end{abstract}

\section*{Introduction}

The parametric degree of a rational surface is the degree of the polynomials
in the smallest possible proper parametrization. 
In the absence of base points, the parametric degree
is just the square root of the degree of the surface.
On the other hand, there are examples
of series of rational surfaces showing that the
number of base points in the smallest parametrization
can be much larger than the degree (see Example~\ref{ex:high}).
It is therefore clear that the degree
does not tell too much about the parametric degree of a rational surface.

Example \ref{ex:torus} below shows that the parametric degree 
is not a geometric but an arithmetic concept, in the sense that it
depends on the choice of the ground field. 
For the complex case, \cite{Schicho:98d} gave upper and
lower estimations for the parametric degree in terms of the degree
of the surface. 

In this paper, we introduce two geometrical invariants of a rational surface,
namely level and keel. These two numbers govern the parametric degree
in the sense that there exist linear upper and lower bounds.
In particular, level and keel determine
the parametric degree up to a multiplicative factor of $2$,
independent on the choice of the ground field (as long as
it is perfect).
We can therefore say that the parametric degree depends ``just slightly''
on the choice of the ground field, as this choice can change by
a factor of at most~2. Of course, this is a bit misleading
because there are surfaces for which rationality depends on the
choice of the field, for instance surfaces that are rational
over $\mathbb{C}$ but not rational over $\mathbb{R}$.

We note that this result implies that the question of rationality
is decidable over any perfect field with decidable first order theory:
first, we compute a parametrization over the algebraic closure,
using \cite{Schicho:97}. If it is of degree $d$, then, by the result 
in this paper, the surface is rational iff there is a parametrization
of degree at most $2d$. For a fixed surface, this question can be 
formulated as a first order sentence. Of course, decidability
of rationality of surfaces is well-known for the real case
by Comesatti's theorem (see \cite{Comesatti:12,Schicho:00c}).

\section{Parametric Degree}

Throughout the paper, we fix field $\mathbb{K}$,
which is assumed to be perfect.
A proper parametri\-zation of a rational surface is a birational map
$\nu:\mathbb{P}^2\to S\subset\mathbb{P}^r$, $r\ge 2$, which is defined
over $\mathbb{K}$. We can write $\nu$ as 
\[ (x_0{:}x_1{:}x_2) \mapsto (F_0(x_0,x_1,x_2){:}\dots{:}F_r(x_0,x_1,x_2)), \]
with $F_0,\dots,F_r\in \mathbb{K}[x_0,x_1,x_2]$ 
homogeneous of the same degree $d$
and without a common divisor. This representation is unique up to the
multiplication by a nonzero constant. The number $d$ is uniquely defined
and it is called the degree of the parametrization. (This should not be
mixed up with the concept of degree of a rational map; in fact, any proper
parametrization is a rational map of degree~1.)

The surface $S$ has, in general, proper parametrizations of different degree.
The smallest possible degree is called the {\em parametric degree}
of $S$, and denoted by $\pdeg(S)$. 
For instance, we have $\pdeg(\mathbb{P}^2)=1$,
and $\pdeg(Q)=2$ when $Q$ is a quadric surface in $\mathbb{P}^2$ with
a $\mathbb{K}$-rational nonsingular point. In the second case, the inverse of
the stereographic projection from this $\mathbb{K}$-rational point defines
a parametrization of degree~2. 

The parametric degree is an arithmetic concept, i.e. it depends
on the choice of the field $\mathbb{K}$. 

\begin{example} \label{ex:torus}
The torus with equation
\[ (x^2+y^2+z^2+\frac{16}{25}w^2)^2-4x^2w^2-4y^2w^2 = 0 \]
has a complex parametrization of degree~3, namely
\[ (s{:}t{:}u) \mapsto 
	((s^2+u^2)(3u-8t):i(s^2-u^2)(3u-8t): \]
\[	-6is(u^2-6tu+8t^2:20st(u-3t)) . \]
The smallest real parametrization has degree~4.
So we have $\pdeg(S_\mathbb{C})=3$ and $\pdeg(S_\mathbb{R})=4$.
\end{example}

\section{Preliminaries}

Most definitions and theorems in this section are well-known.
The single exception is
the definition of the function $\mathrm{nmc}$ at the end of the section,
which turns out to be convenient for the definition of the keel.

Let $S$ be a nonsingular projective surface over $\mathbb{K}$.
\begin{itemize}
\item A prime divisor is an irreducible curve on $S$. Note that it
	is required that the curve is defined over $\mathbb{K}$,
	but it may split into several components over the algebraic
	closure.
\item $\mathrm{Div}(S)$, the group of divisors of $S$, is the free
	abelian group generated by the prime divisor. Divisors are
	denoted by captital letters.
\item A divisor is called effective iff all its coefficients
	are nonnegative. If $B-A$ is effective, then we also
	say that $A$ divides $B$ or $A\le B$. 
\item Two effective divisor without common component are equivalent
	iff they are two fibers of a rational map $S\to\mathbb{P}^1$.
	Linear equivalence is the finest equivalence relation 
	on $\mathrm{Div}(S)$ which is compatible with addition
	and for which the previous statement is true.
\item $\mathrm{Pic}(S)$ is the group of classes of divisors.
	Classes are denoted by capital letters.
\item The class $D$ is called effective iff it has an effective divisor.
\item $|D|$ is the set of all effective divisors in the class $D$.
	It has a natural structure of a projective space over $\mathbb{K}$.
	The dimension of this projective space is denoted by $\dim(D)$.
	A linear system of divisors is a subset of $|D|$ corresponding
	to a projective subspace.
\item If the point $p\in{S}$ is contained in all divisors
	in a linear system $l$, then $p$ is a base point of $l$. 
	If the complete linear system $|D|$ has no
	base points, then we say that $D$ is free.
\item If $l$ is a non-empty linear system,
	then the associated rational map is denoted by
	$\phi_l:S\to\mathbb{P}^{\dim(l)}$. It is defined outside
	the base locus. The codomain $\mathbb{P}^{\dim(l)}$ is naturally
	identified with the dual projective space of $l$; the image
	of $p$ corresponds to the subset of divisors in $l$ passing
	through $p$, which is a hyperplane.
\item The intersection product $\mathrm{Pic}(S)^2\to\mathbb{Z}$
	is symmetric and bi-additive, and if the classes $A$ resp.
	$B$ contain two effective divisors $A_0$ resps. $B_0$ without
	common component, then $AB$ is the number of common points
	of $A_0$ and $B_0$, properly counted. Especially, $AB\ge 0$
	in this case.
\item $K$ or $K_S$ is the canonical class of $S$.
\item A class $D$ is called nef iff $DC\ge 0$ for all effective $C$.
\end{itemize}

We also need the following well-known theorems.

\begin{theorem} \label{thm:c1}
A birational regular map $\phi:S_1\to S_2$ induces two homomorphisms,
the pushforward $\phi_\ast:\mathrm{Div}(S_1)\to\mathrm{Div}(S_2)$ and
the pullback $\phi^\ast:\mathrm{Div}(S_2)\to\mathrm{Div}(S_1)$.
Both functions are well defined on classes and preserve effectivity. 

For $C\in\mathrm{Pic}(S_1)$ and $D\in\mathrm{Pic}(S_2)$, the following hold.
\begin{itemize}
\item $(\phi_\ast\circ\phi^\ast)(D)=D$. 
\item $\dim(\phi_\ast(C))\ge\dim(C)$. 
\item $\dim(\phi^\ast(D))=\dim(D)$. 
\item $\phi_\ast(C)D=C\phi^\ast(D)$. 
\item $(\phi_\ast(C))^2\ge C^2$. 
\item $(\phi^\ast(D))^2= D^2$. 
\item $\phi_\ast(K_{S_1})=K_{S_2}$.
\end{itemize}
\end{theorem}

\begin{theorem} \label{thm:c2}
Let $S$ be a nonsingular projective surface. Let $E\in\mathrm{Div}(S)$
be a prime divisor such that $E^2=EK<0$. A prime divisor with these
properties is called exceptional divisor.

Then there exists
a regular birational map $\pi:S\to S'$, called the blowing down of $E$,
such that the kernel of $\pi_\ast$ is generated by $E$. 
Moreover, $K_{S}=\pi^\ast(K_{S'})+E$.

Any birational regular map is a composition of such blowing down maps.
\end{theorem}

We define a function $\mathrm{nmc}$ 
(for ``number of moving components'')
from effective classes to nonnegative integers.
Let $X(D)\subset\mathbb{P}^{\dim(D)}$ be the image of the associated
rational map $\phi_D$.
\begin{itemize}
\item If $X(D)$ is a point (this is the case iff $\dim(D)=0$),
	then $\mathrm{nmc}(D):=0$.
\item If $X(D)$ is a curve of degree $m$, then $\mathrm{nmc}(D):=m$.
\item If $X(D)$ is a surface, then $\mathrm{nmc}(D):=1$.
\end{itemize}

\section{Level and Keel}

We first introduce level and keel for divisors on nonsingular
surfaces. Then the concepts are transferred to embedded surfaces
with arbitrary singularities, using a resolution of singularities.

Let $D$ be an effective divisor class of $S$.
The {\em level} of $D$ is the supremum of all rational numbers
$p/q$, $q>0$, such that $qD+pK$ is effective.
If the supremum is assumed, then the {\em keel} of $D$ is equal to
the supremum of all numbers of the form $\frac{\mathrm{nmc}(qD+pK)}{q}$
where $p/q$ is the level.
If the supremum in the definition of the level is not assumed 
(for instance if the level is irrational or infinity), 
then the keel is defined as $0$.

\begin{remark}
If some multiple of $K$ is effective, then we have $\mathrm{level}(D)=\infty$
and $\mathrm{keel}(D)=0$. Hence level and keel are only interesting
if the Kodaira dimension of $S$ is negative.
\end{remark}

\begin{remark}
The numbers $\dim(D)$ and $\mathrm{nmc}(D)$ are preserved under extension
of the ground field. It follows that level and keel are preserved
under extension of the ground field (i.e. they are geometric).
\end{remark}

\begin{example} \label{ex:p2d}
Let $S=\mathbb{P}^2$ and $D=nL$, where $L$ is the class of lines and $n\ge 0$.
Then $K=-3L$, and a class $mL$ is effective iff $m\ge 0$, and we
have $\mathrm{nmc}(0)=0$.
It follows that $\mathrm{level}(D)=n/3$ and $\mathrm{keel}(D)=0$.
\end{example}

\begin{example} \label{ex:pp}
Let $S=\mathbb{P}^1\times\mathbb{P}^1$ and $D=mF_1+nF_2$, where $F_1,F_2$ 
are the classes of the fibers of the two projections and $0\le m\le n$.
Then $K=-2F_1-2F_2$, and a class $aF_1+bF_2$ is effective iff 
$a\ge 0$ and $b\ge 0$, and we have $\mathrm{nmc}(aF_2)=a$.
It follows that $\mathrm{level}(D)=m/2$ and $\mathrm{keel}(D)=n-m$.
\end{example}

\begin{example} \label{ex:tor}
We sketch an example which generalizes the two above.

\begin{figure}[htbp]
  \centering
  \includegraphics{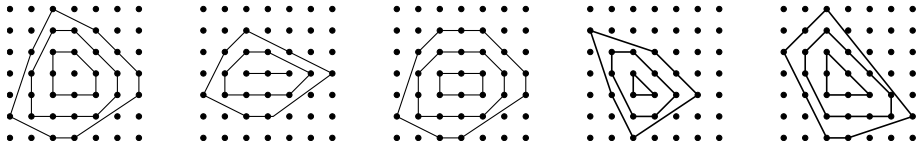}
  \caption{Polygons of level and keel (3,0), 
	(2,2), (5/2,1), (7/3,0), and (8/3,0).}
  \label{fig:latt}
\end{figure}

Let $\Gamma$ be a convex lattice polygon, i.e.
the convex hull of a finite number of points in the plane
with integer coordinates. 
The polygon $\Gamma$ defines a nonsingular toric surface $S$
(the minimal resolution of the toric surface defined by the inner normals)
and an effective divisor
(the inverse image of the class of hyperplane section in the
projective embedding defined by $\Gamma$).

The class $qD+pK$ corresponds to the convex figure obtained by
scaling $\Gamma$ by a factor of $q$ and moving each edge $p$ steps
inward. The class is effective iff this figure is non-empty.
Hence the level is equal to $p/q$ if we can enlarge $\Gamma$ by
a factor of $q$, pass $p$ times to the convex hull of the interior
points, and obtain a line segment or a point 
(see figure~\ref{fig:latt}; cf also \cite{Schicho:04b}).
The keel is the number of points on this line segment or point, minus 1,
divided by $q$.
\end{example}

We now define level and keel of a projective surface $S\subset\mathbb{P}^r$
(possibly with singularities). Let $\pi:\tilde{S}\to S$
be a resolution of singularities, i.e. a proper birational map
such that $\tilde{S}$ is nonsingular. 
Let $H\in\mathrm{Pic}(\tilde{S})$ be the class of the pullbacks
of hyperplane sections.
We define the level and the keel of $S$ as the level and the keel
of $H$. Theorem~\ref{thm:wd} below says that this
is independent of the choice of the desingularization.

\begin{lemma} \label{lem:fp}
Let $\phi:S_1\to S_2$ be a regular birational map.
Let $D\in\mathrm{Pic}(S_2)$.
Let $C$ be an effective divisor of $S_1$ such that $\phi_\ast(C)=0$.
Then $C$ is a common divisor of the linear system $|\phi^\ast(D)+C|$.
\end{lemma}

\begin{proof}
The proof proceeds by induction on the number of blowing downs
into which the birational regular map $\phi$ can be decomposed.
First, assume that $\phi:S_1\to S_2$ is the blowing down
of an exceptional divisor $E$. Then $C=nE$ for some $n\ge 0$.
Let $m\ge 0$ be the largest number such that $mE$ is a common divisor of
$|\phi^\ast(D)+nE|$. Then there exists an effective divisor 
in $|\phi^\ast(D)+(n-m)E|$ that does not have $E$ as component.
It follows that $(\phi^\ast(D)+(n-m)E)E=m-n\ge 0$, which shows
that $nE$ is a common divisor.

Second, assume that $\phi$ can be decomposed as
$S_1\stackrel{\phi_1}{\to}S_3\stackrel{\phi_2}{\to}S_2$,
where $\phi_1$ is the blowing down of an exceptional divisor $E$.
Because the pushforward preserves effectivity, ${\phi_1}_\ast(C)$
is effective. By induction, ${\phi_1}_\ast(C)$ is
a common divisor of the linear system $|{\phi_2}^\ast(D)+{\phi_1}_\ast(C)|$. 
Because the dimension of a linear system is preserved by pullback,
the equation 
$\dim({\phi_2}^\ast(D)+{\phi_1}_\ast(C))=\dim({\phi_2}^\ast(D))$
implies the equation
$\dim({\phi}^\ast(D)+{\phi_1}^\ast{\phi_1}_\ast(C))
	=\dim({\phi}^\ast(D))$,
which implies that ${\phi_1}^\ast{\phi_1}_\ast(C)$ is a common divisor of
$|{\phi}^\ast(D)+{\phi_1}^\ast{\phi_1}_\ast(C)|$. 
Now $C-{\phi_1}^\ast{\phi_1}_\ast(C)$ lies in the kernel of ${\phi_1}_\ast$.
Therefore it is a multiple $nE$ of the exceptional divisor.

We distinguish two cases. If $n\le 0$, 
then $C\le{\phi_1}^\ast{\phi_1}_\ast(C)$, and it follows
\[ \dim({\phi}^\ast(D)) \le \dim({\phi}^\ast(D)+C) \le 
	\dim({\phi}^\ast(D)+{\phi_1}^\ast{\phi_1}_\ast(C)) =
	\dim({\phi}^\ast(D)) , \]
hence we have equality everywhere and the statement is proved.

If $n>0$, then $nE$ is effective, and $nE$ is a common divisor of the linear
system $|{\phi_1}^\ast({\phi_2}^\ast(D)+{\phi_1}_\ast(C))+nE|$
by the induction base case. Therefore, we have
\[ \dim({\phi}^\ast(D)+C) =
	\dim({\phi_1}^\ast({\phi_2}^\ast(D)+{\phi_1}_\ast(C))+nE)= \]
\[ \dim({\phi_1}^\ast({\phi_2}^\ast(D)+{\phi_1}_\ast(C)))=
	\dim({\phi}^\ast(D)) , \]
hence the statement is also proved.
\end{proof}

\begin{theorem} \label{thm:wd}
Let $\pi_1:\tilde{S}_1\to S$ and $\pi_2:\tilde{S}_2\to S$ be two
desingularizations of $S$. Let $H_1\in\mathrm{Pic}(\tilde{S}_1)$
and $H_2\in\mathrm{Pic}(\tilde{S}_2)$ be the pullbacks of hyperplane
sections. Then we have 
\[ \mathrm{level}(H_1)=\mathrm{level}(H_2), \ 
	\mathrm{keel}(H_1)=\mathrm{keel}(H_2). \]
\end{theorem}

\begin{proof}
First, let us assume that there exists a birational regular map
$\phi:\tilde{S}_1\to\tilde{S}_2$ such that $\pi_1=\pi_2\circ\phi$.
Then $\phi$ transforms hyperplane pullbacks to hyperplane pullbacks,
i.e. $\phi^\ast(H_1)=H_2$.
Let $p,q$ be positive integers.
By Theorem~\ref{thm:c2}, the class $C:=K_1-\phi^\ast(K_2)$ is effective.
By Lemma~\ref{lem:fp}, the divisor $pC$ is a common divisor of the linear
system $|\phi^\ast(qH_2+pK_2)+pC|=|qH_1+pK_1|$.
It follows that the two linear systems $|qH_2+pK_2|$ and
$|qH_1+pK_1|$ have the same dimension and the same number of
moving components, and the statement is proven.

In the general case, there exists a dominating desingularization
$\pi_3:\tilde{S}_3\to S$ and birational regular maps
$\phi_i:\tilde{S}_3\to\tilde{S}_i$ such that $\pi_3=\pi_i\circ\phi_i$
for $i=1,2$. Hence it can be reduced to the special case above.
\end{proof}

\begin{example} \label{ex:p2}
Assume that $\nu:\mathbb{P}^2\to S$ is a parametrization of degree $d$
without base points. Then $\nu$ is regular and we can use it
as resolution of singularities. It follows that $H=dL$ and
we have $\mathrm{level}(S)=\mathrm{level}(H)=d/3$ and
$\mathrm{keel}(S)=\mathrm{keel}(H)=0$, by Example~\ref{ex:p2d}.
\end{example}

It is also convenient to extend the notion of parametric degree to
divisors. Let $S$ be a nonsingular surface. Let $D\in\mathrm{Pic}(S)$
be an nef divisor. 
A linear system $l$ of divisors
is called parametrizing iff $\dim(l)=2$ and 
$\phi_l:S\to\mathbb{P}^2$ is birational.
A class $P$ is called parametrizing iff
$|P|$ contains a parametrizing linear system.
Then we define $\pdeg(D)$ as the minimum
of all numbers $PD$, where $P$ is a parametrizing class.

\begin{lemma} \label{lem:od}
Let $S\subset\mathbb{P}^r$ be a (possibly singular) surface.
Let $\pi:\tilde{S}\to S$ be a resolution of its singularities.
Let $H\in\mathrm{Pic}(\tilde{S})$ be the class of the pullbacks
of hyperplane sections.
Then $\pdeg(S)=\pdeg(H)$.
\end{lemma}

\begin{proof}
There is a one-to-one correspondence of parametrizations of $S$
and pa\-ra\-met\-ri\-zing linear systems of $D$, and the degrees coincide
for corresponding parametrization/class.
\end{proof}

\begin{remark} \label{rem}
Lemma~\ref{lem:od} allows to reduce any relation between parametric
degree, level and keel of a singular surfaces to the same relation
between parametric degree, level and keel of a divisor on a nonsingular
surface.
\end{remark}

\section{The Lower Bound}

The main idea for establishing a lower bound for the parametric
degree in terms of level and keel is to analyze what happens in
the examples~\ref{ex:p2d} and \ref{ex:p2} when the parametrization
has base points.

\begin{theorem} \label{thm:lb}
Let $S\subset\mathbb{P}^r$ be a rational surface. Then we have
\[ \pdeg(S) \ge 3\ \mathrm{level}(S)+\mathrm{keel}(S) . \]
\end{theorem}

\begin{proof}
It suffices to prove the above inequality for a divisor $H$.
We assume that $PH=\pdeg(H)$ and that some linear subsystem of $|P|$
induces a birational map to $\mathbb{P}^2$. Moreover,
we assume that $l\subseteq |P|$ is a parametrizing linear system.
If $C$ is a common divisor of $l$, then we would have $(P-C)H\le PH$
because $H$ is nef. In this case, we can replace $P$ by $P-C$
which is also a parametrizing class. As we can do this only finitely
many times, because $l$ has only finitely many common components,
we can assume that $l$ has no common components.

Claim 1: $P$ is nef. Indeed, if $C$ is a prime divisor, then
$CP\ge 0$ because $C$ is not a common component of $l$.

Claim 2: if $C$ is a prime divisor with positive dimension,
then $PC>0$. Indeed, we have $P^2>0$ because the image of $\phi_l$
is $\mathbb{P}^2$, and therefore two generic divisors in $l$
intersect in a point outside the base locus. And $|C|$ contains
two divisors without common component, hence $C^2\ge 0$. Therefore $PC$
cannot be zero by the Hodge index theorem.

Claim 3: $PK\le -3$. To prove this, we resolve the base points
of $\phi_l:S\to\mathbb{P}^2$ and get a birational regular map
$\pi:\tilde{S}\to S$ such that $\phi_P\circ\pi:\tilde{S}\to\mathbb{P}^2$
is regular. Then $\pi^\ast(P)$ is nef by Theorem~\ref{thm:c2}, and it follows
\[ PK = \pi^\ast(P)\pi^\ast(K)=\pi^\ast(P)K_{\tilde{S}}\le
	\pi^\ast(P)(\phi_l\circ\pi)^\ast(K_{\mathbb{P}^2}) \]
\[ = (\phi_l\circ\pi)_\ast(\pi^\ast(P))K_{\mathbb{P}^2} 
	= L(-3L)=-3 , \]
where $L\in\mathrm{Pic}(\mathbb{P}^2)$ is the class of lines.

For two positive integers $p,q$, the divisor $qH+pK$ can be effective
only if $(qH+pK)P=q\pdeg(H)-3p\ge 0$. This proves that 
$\mathrm{level}(H)\le \pdeg(H)/3$.

Now, assume $p/q=\mathrm{level}(H)$. We claim that 
$\mathrm{nmc}(qH+pK)\le \pdeg(H)q-3p$. Let $F$ be the greatest
common divisor of $qH+pK$, and let $B$ be a generic divisor
in $|qH+pK-F|$. Then $B$ corresponds
to a generic hyperplane section of the associated image $X:=\phi_{qH+pK}(S)$.
If the image $X$ is a point, then $\mathrm{nmc}(qH+pK)=0$ and
the claim is true. If $X$ is a surface, then $B$ has positive
dimension; it follows that 
$\pdeg(H)q-3p\ge BP\ge 1$, and because $\mathrm{nmc}(qH+pK)=1$
the claim is true. If $X$ is a curve, then $X$ is necessarily rational,
and $\phi_{qH+pK}$ factors through a rational map
$S\to\mathbb{P}^1$ which is associated to some divisor $A$.
Because the divisors in $A$ are also the fibers of 
$\phi_{qH+pK}$, we have $B=mA$, where $m=\mathrm{nmc}(qH+pK)$ is the
number of intersection points of $X$ with a generic hyperplane.
Since $A$ has positive dimension, we have 
\[ \pdeg(H)q-3p\ge BP = mAP \ge m=\mathrm{nmc}(qH+pK) ,\]
hence the claim is true also in this last case.
This shows that 
\[ \mathrm{keel}(H)\le \pdeg(H)-3p/q 
	= \pdeg(H)-\mathrm{level}(H) . \]
\end{proof}

\begin{remark}
By analyzing Example~\ref{ex:pp} in a similar fashion, 
one can show that if $S$ has a parametrization
of bidegree $(m,n)$, $m\le n$, then $m\ge 2\ \mathrm{level}(S)$
and $n\ge 2\ \mathrm{level}(S)+\mathrm{keel}(S)$. 
\end{remark}

\begin{example} \label{ex:high}
Here is an example that shows how to use the concepts
of level and keel in order to construct surfaces that have
a high parametric degree.

Let $n\ge 5$ be an odd integer. Let $S\subset\mathbb{P}^3$
be the surface given by the equation
\[ z^{2n+1}-x^2w^{2n-1}-y^nw^{n+1} = 0 . \]
We can compute the level and keel by using the parametrization
\[ (s{:}t{:}u) \mapsto 
	((s^2u^{n-2}+t^n)^ns:(s^2u^{n-2}+t^n)^2tu^{n^2-2n} \]
\[ :(s^2u^{n-2}+t^n)u^{n^2-n+1}:u^{n^2+1}) . \]
By resolving the base points of the parametrization, we get a 
resolution $\tilde{S}$ of the singularities of $S$. Explicit analysis
of the base points (see \cite{Schicho:02a}) shows that there is
one base point of multiplicity $n^2-2n$, with $\frac{n-3}2$ base points
of multiplicity $2n$ and $2n+3$ base points of multiplicity $n$ in the
infinitely near, and one simple base point with $n^2-2n-1$ simple
base points in the infinitely near.

For positive integers $p,q$, the linear system $|qH+pK|$ on $\tilde{S}$
corresponds to the linear space of forms of degree
$q(n^2+1)-3p$ that vanish with multiplicity $qr-p$ at each point
of multiplicity $r$. Such forms exist for $2p\le (2n+1)q$,
hence $\mathrm{level}(S)=n+\frac{1}{2}$. If $2p=(2n+1)q$, then
the corresponding linear space is the vectorspace of forms
of degree $\frac{2n^2-6n-1}{2}q$ vanishing with multiplicity
$\frac{2n^2-6n-1}{2}q$ at the $(n^2-2n)$-fold base point and
with multiplicity $\frac{2n-1}{2}q$ at the $\frac{n-3}2$ base points
of multiplicity $2n$ in the infinitely near. Hence
\[ \mathrm{keel}(S)=\frac{2n^2-6n-1}{2}-\frac{n-3}{2}\ \frac{2n-1}{2}
	=\frac{2n^2-5n-5}4 . \]

By Theorem~\ref{thm:lb}, the parametric degree is greater than
or equal to $\frac{2n^2+7n+1}{4}$. Because there is a parametrization
of degree $n^2+1$, we know that the parametric degree grows
proportional to the square of the implicit degree.
\end{example}

\section{The Adjoint Chain}

Adjoints are a tool for constructing minimal models of a given
surface or higher-dimensional varieties. Starting with a nef
divisor class, we keep alternating to blow down orthogonal exceptional
divisors and adding the canonical class, until nefness does
not hold any more. The last surface with nef class in the
process is a minimal model with special properties. This technique
has been used in \cite{Schicho:97} to construct parametrizations
in the case $\mathbb{K}$ is algebraically closed. Similar 
constructions appear in various other contexts,
see \cite{Beltrametti_Sommese:95} for a survey.

Let $S$ be a nonsingular surface, and let $D\in\mathrm{Pic}(S)$.
Following \cite{Manin:66,Manin:67}, 
we say that $S$ is $D$-minimal iff $S$ has no exceptional divisor
orthogonal to $D$. We say $S$ is minimal iff $S$ is $0$-minimal.
The following theorem is well-known (see \cite{Manin:66}).

\begin{theorem}
Let $S$ be a nonsingular surface and let $D\in\mathrm{Pic}(S)$.
Then there exists a birational regular map $\mu:S\to S_0$,
such that $D=(\mu^\ast\circ\mu_\ast)(D)$ and $S_0$ is
$\mu_\ast(D)$-minimal. We call this a $D$-minimalization.
\end{theorem}

In the rest of the paper, we fix a nonsingular
rational surface $S$ and
a nef class $D\in\mathrm{Pic}(S)$ such that $D^2>0$ (for instance,
the class of pullback of hyperplane sections in a resolution).
The adjoint chain $\mathcal{S}$ 
is a chain of surfaces and birational regular maps
\[ S\stackrel{\mu_0}{\to} S_0 \stackrel{\mu_1}{\to} 
	S_1 \stackrel{\mu_2}{\to} \dots \]
and divisor classes $D_i\in\mathrm{Pic}(S_i)$, which is constructed recursively
in the following way. First, we let $\mu_0:S\to S_0$ be a $D$-minimalization
of $S$, and we let $D_0:={\mu_0}_\ast(D_0)$. 
Now assume that we have already defined $S_i$ and $D_i$.
Let $K_i$ be the canonical class of $S_i$.
If $D_i+K_i$ is not effective, then the adjoint chain ends;
we denote the index of the last surface with $a$.
Otherwise, we let $\mu_{i+1}:S_i\to S_{i+1}$ be a $(D_i+K_i)$-minimalization 
of $S_i$, and we let 
$D_{i+1}:=(\mu_{i+1})_\ast(D_i+K_i)$. If the adjoint chain is
infinite, then we set $a:=\infty$ (but we will prove that $a$ is finite).

\begin{lemma}
The classes $D_i$ above are effective and nef. 
If $i<a$, then ${D_i}^2>0$.
\end{lemma}

\begin{proof}
If $\mathbb{K}$ is algebraically closed, then the proof is
well-known (\cite{Schicho:97}, Lemma~A.2 and Lemma~A.3).
There is only one step in the proofs of \cite{Schicho:97}
that uses the assumption that $\mathbb{K}$ is algebraically closed:
in this case, any prime divisor $C$ of dimension~0 with $CK<0$
is exceptional, and we have $C^2=CK=-1$
(Lemma~A.1 in \cite{Schicho:97}). It follows that if $D$ is nef
and $D+K$ is effective but not nef, then there exists an exceptional divisor
orthogonal to $D$.

Here is an adaption of the proof to the case of non-closed fields:
assume that $D$ is nef and $D+K$ is effective but not nef. 
Let $C$ be a prime divisor such that $(D+K)C<0$. 
Then $CK<0$, and $\dim(C)=0$ because $C$ must be fixed in $|D+K|$.
By the lemma below, $C$ is exceptional and $CD$ is an integral
multiple of $CK$. This is only possible if $CD=0$, hence $C$ is
orthogonal to $D$, and the rest of the proof works as in the case
where $\mathbb{K}$ is algebraically closed.
\end{proof}

\begin{lemma}
Let $C$ be a prime divisor such that $\dim(C)=0$ and $CK<0$. Then $C$ is 
exceptional, and $CD$ is an integral multiple of $CK$ 
for all $D\in\mathrm{Pic}(D)$.
\end{lemma}

\begin{proof}
Let $\overline{S}$ be the surface obtained by base field extension
to the algebraic closure $\overline{\mathbb{K}}$. There are natural
injections $\mathrm{Div}(S)\to \mathrm{Div}(\overline{S})$
and $\mathrm{Pic}(S)\to \mathrm{Pic}(\overline{S})$.
In general, $C$ need not be a prime divisor in $\mathrm{Div}(\overline{S})$,
but it has only simple components $C=\sum_{i=1}^r C_i$. Each $C_i$
has dimension~0. Moreover, $C_iD=C_jD$ for
any $i,j\le r$ and $D$ coming from $S$, 
because $C_i$ and $C_j$ are conjugate under the action
of the Galois group of the extension $\mathbb{K}\subset\overline{\mathbb{K}}$. 
Especially, $C_iK=C_jK$.
It follows that $C_iK<0$.
By Lemma~A.1 in \cite{Schicho:97}, each $C_i$ is exceptional,
and $C_i^2=C_iK=-1$ for all $i$. Hence 
$CK=-r$, and $CD$ is
an integral multiple of $r$.

It remains to show that $C^2=-r$. For $i\ne j$, we have 
$\dim(C_i+C_j)\ge C_iC_j$ by the Riemann-Roch theorem. On the
other hand, $\dim(C_i+C_j)=0$ because $\dim(C)=0$,
hence $C_iC_j=0$. Hence $C^2=\sum_{i=1}^rC_i^2=-r$.
\end{proof}

\begin{lemma} \label{lem:pqi}
Let $p,q,i$ be integers such that $i\le a$ and $p\ge qi$.
Then the linear systems $|qD+pK|$ and $|qD_i+(p-qi)K_i|$ have the
same dimension and number of moving components.
\end{lemma}

\begin{proof} 
For $j=0,\dots,i$, 
let $\phi_j:S\to S_j$ be the map $\mu_j\circ\dots\circ\mu_0$.
Then $\phi_j^\ast(D_j)=\phi_{j-1}^\ast(D_{j-1}+K_{j-1})$.
The class $E_j:=\phi_j^{\ast}K_j-\phi_{j-1}^{\ast}K_{j-1}$ is effective.
Therefore, the class 
\[ C:=qD+pK-\phi_i^\ast(qD_i+(p-qi)K_i)=
        \sum_{j=0}^i (p-qj)E_j \]
is effective, too. By Lemma~\ref{lem:fp}, the unique divisor in $|C|$
is fixed in $|qD+pK|$. Because
the pullback preserves the dimension and the number of moving
components, the lemma follows.
\end{proof}

\begin{corollary} \label{cor}
We have $\mathrm{level}(D_i)=\mathrm{level}(D)-i$
and $\mathrm{keel}(D_i)=\mathrm{keel}(D)$.
\end{corollary}

\begin{lemma}
We have $a\le\mathrm{level}(D)$. In particular, $a$ is finite.
\end{lemma}

\begin{proof}
By Corollary~\ref{cor}, it suffices to prove that $\mathrm{level}(D_a)\ge 0$.
But this is clear since $D_a$ is effective.
\end{proof}

The following lemmas can be used to compute level and keel in terms
of the adjoint chain.

\begin{lemma} \label{lem:c0}
Assume that $D_a=0$.
Then $\mathrm{level}(D)=a$ and $\mathrm{keel}(D)=0$.
\end{lemma}

\begin{proof}
This follows immediately from Corollary~\ref{cor} and from
$\mathrm{level}(0)=0$ and $\mathrm{keel}(0)=0$.
\end{proof}

\begin{lemma} \label{lem:c1}
Assume that ${D_a}^2=0$ and $D_a\ne 0$.
Then $\mathrm{level}(D)=a$ and $\mathrm{keel}(D)=\mathrm{nmc}(D_a)=:k>0$,
and there is a free divisor $P\in\mathrm{Pic}(S_a)$ such that
$D_a=kP$ and $PK_a=-2$ and $P^2=0$ and $\dim(P)=1$.
\end{lemma}

\begin{proof}
The proof of Lemma~A.7 in \cite{Schicho:97} generalizes
without problems to non-closed fields. This shows the existence
of $P$ with the desired properties.

If $p,q>0$, then $qD_a+pK_a$ cannot be effective because
$(qD_a+pK_a)P=-2p<0$. Hence $\mathrm{level}(D_a)=0$ and
$\mathrm{level}(D)=a$ by Corollary~\ref{cor}. Moreover,
$\mathrm{nmc}(qD_a)=\mathrm{nmc}(qkP)=qk$, hence 
$\mathrm{keel}(D_a)=\mathrm{keel}(D)=k$.
\end{proof}

\begin{lemma} \label{lem:c2}
Assume that ${D_a}^2>0$. 
Then one of the following cases holds. \\
a) $\mathrm{level}(D)=a+1/3$, $\mathrm{keel}(D)=0$, and $3D_a+K_a=0$. \\
b) $\mathrm{level}(D)=a+2/3$, $\mathrm{keel}(D)=0$, and $3D_a+2K_a=0$. \\
c) $\mathrm{level}(D)=a+1/2$, $\mathrm{keel}(D)=0$, and $2D_a+K_a=0$. \\
d) $\mathrm{level}(D)=a+1/2$, $\mathrm{keel}(D)=\mathrm{nmc}(2D_a+K_a)/2>0$,
	and $(2D_a+K_a)^2=0$. \\
In particular, level and keel are rational numbers with a denominator
dividing~6.
\end{lemma}

\begin{proof}
Lemma~A.8 from \cite{Schicho:97} -- which is also true in the case
$\mathbb{K}$ is non-closed -- says that we can conclude from
${D_a}^2>0$ that either $3D_a+K_a=0$, or $3D_a+2K_a=0$,
or $(2D_a+K_a)^2=0$. Using Corollary~\ref{cor}, we get
$\mathrm{level}(D)=\mathrm{level}(D_a)+a=a+i/3$ in the $i$-th case,
for $i=1,2$. In the third case, we apply Lemma~\ref{lem:c0} or
Lemma~\ref{lem:c1} to $D:=2D_a+K_a$, and we get either (c) or (d), 
depending whether $2D_a+K_a$ is zero or not.
\end{proof}

\begin{remark}
At this point, it is instructive to revisit Example~\ref{ex:tor}
again. Starting from a convex lattice polygon, we pass to the
convex hull of the interior points $a$ times. If we obtain a single
point, then $D_a=0$ holds for the corresponding toric surface.
If we obtain a line segment with $k+1$ lattice points, then $D_a=kP$
for some $P$ with $P^2=0$ and $PK_a=-2$, as in Lemma~\ref{lem:c1}.
If we obtain a lattice polygon without interior lattice points, 
then one of the following four cases holds: \\
a) after scaling by 3, we get a polygon with one interior point; \\
b) after scaling by 3 and passing to the convex hull of interior points, 
we get a polygon with one interior point; \\
c) after scaling by 2, we get a polygon with one interior point; \\
d) after scaling by 2, we get a polygon with several interior points
that are all on a line. \\
Of course, these are instances of the 4 cases (a), (b), (c), (d)
in Lemma~\ref{lem:c2}.
\end{remark}

\begin{remark}
The lemmas above remind on the Kawamata Rationality Theorem 
and the Kawamata-Shokurov Base Point Free Theorem 
(see \cite{Kawamata_Matsuda_Matsuki:85}):
if $D$ is ample, then the ``nefness value'' $v$ is a rational number,
and some multiple of $D+vK$ is free. The associated contraction
morphism is either a blowing down, or a map with conic fibers,
or a constant map. Of course, the Kawamata Rationality Theorem
and the Kawamata-Shokurov Base Point Free Theorem
hold in a much more general context (arbitrary dimension,
rationality need not be assumed).
\end{remark}

\begin{remark}
If $S$ is a rational surface with degree $d$ and sectional genus $p_1$,
then we have the inequality $a+\dim(D_a)\le p_1+{2p_1-d-1 \choose 2}$,
by Lemma~8 in \cite{Schicho:98d}. Using the classification
of surfaces occuring in Lemma~\ref{lem:c2} 
(see Lemma~A.8 in \cite{Schicho:97}), it is easy to check that 
$\mathrm{keel}(D)\le \dim(D_a)$ in all cases. Together with the upper bound
$p_1\le {d-1 \choose 2}$ for the sectional genus, we get the bound
$\mathrm{level}(D)+\mathrm{keel}(D)\le d^4/2$. Together with
the bound in Theorem~\ref{thm:fin} below, we obtain the bound
$\pdeg(S)\le 3\deg(S)^4$.
\end{remark}

\section{The Upper Bound}

In order to establish an upper bound for the parametric degree,
one has to construct a parametrization (or, equivalently,
a parametrizing divisor class). The idea is to construct
a minimal model using adjoints, and then to use the 
well-known classification of such minimal surfaces,
due to Manin and Iskovskih~\cite{Manin:66,Manin:67,Iskovskih:67,%
Iskovskih:70,Iskovskih:72,Iskovskih:80}.

\begin{theorem} \label{thm:i1}
Let $S$ be a minimal rational surface such that $-K$ is nef and $K^2>0$.
Then one of the following cases holds. \\
a) $S\cong\mathbb{P}^2$; in this case, $K=-3L$, where $L$
	is the class of lines. \\
b) $S$ is isomorphic to a quadric in $\mathbb{P}^3$ or to the blowup
	of a singular quadric cone in $\mathbb{P}^3$ at its vertex. 
	If $Q$ is the
	class of conic plane sections, the $K=-2Q$. \\
c) $S$ is isomorphic to a Del Pezzo surface of degree~5 in $\mathbb{P}^5$.
	Its Picard group is cyclic, generated by $K$.
	The class of hyperplane sections is $-K$. \\
d) $S$ is isomorphic to a Del Pezzo surface of degree~6 in $\mathbb{P}^6$.
	Its Picard group is again cyclic, generated by $K$,
	and the class of hyperplane sections is $-K$.
\end{theorem}

\begin{proof}
Let $d:=K^2$. By the classification of Del Pezzo surfaces over
algebraically closed field (see \cite{Manin:74}, Theorem~24.4), 
we have $1\le d\le 9$.

If $d=9$, then $S$ is a Severi-Brauer surface. As $S$ also has
a parametrization over $\mathbb{K}$, it is isomorphic to $\mathbb{P}^2$
and (a) holds.

If $d=8$, then $S$ is isomorphic to a ruled surface $F_n$ over
the algebraic closure $\overline{\mathbb{K}}$, where $0\le n\le 2$. 
The case $n=1$ is not possible, because in this case $S$ would not
be minimal; in the two remaining cases we have $-K=2Q$ for some
divisor $Q$, whose associated image is a quadric in $\mathbb{P}^3$.

If $d=7$, then $S$ is not minimal.

If $d=6$, then (d) holds.

If $d=5$, then (c) holds.

If $d=4$, then $S$ cannot be both minimal and rational over $\mathbb{K}$,
by Theorem~1 in \cite{Iskovskih:72}.

If $d=3$, $2$, or $1$, then $S$ cannot be both minimal and 
rational over $\mathbb{K}$,
by Theorem~5.7 in \cite{Manin:67}.
\end{proof}

\begin{theorem} \label{thm:j1}
Let $S$ be a nonsingular rational surface and let $D$ be a nef divisor 
such that $\mathrm{keel}(S)=0$ and $D^2>0$. Then
\[ \pdeg(D) \le 6\ \mathrm{level}(D) . \]
\end{theorem}

\begin{proof}
By lemmas~\ref{lem:c0}, \ref{lem:c1}, \ref{lem:c2}, we can reduce
to the case $\mathrm{level}(D)=a$ and $D_a=0$ by replacing $D$
by $2D$ or $3D$. Then $-K_a$ is nef and $(-K_a)^2>0$ because 
$-K_a$ is the direct image of $D_{a-1}$. Then $S_a$ satisfies
the assumptions in Theorem~\ref{thm:i1}.
For each of the cases (a), (b), (c), we construct below a parametrizing class 
$P\in\mathrm{Pic}(S_a)$, such that $P(-K_a)\le 6$. 
Let $\phi_a:S\to S_a$ be the minimalization map.
Then $\phi_a^\ast(P)$ is a parametrizing divisor for $S$, and
\[ \phi_a^\ast(P)D = P{\phi_a}_\ast(D) = P(D_a-aK_a) \le 6a. \]

Case (a): we take $P:=L$. Then $P(-K_a)=3\le 6$.

Case (b): we take $P:=Q$. This is a parametrizing class because
we can choose a parametrizing system $l$ as the linear system 
of conic sections through
a fixed nonsingular point $p$ defined over $\mathbb{K}$. 
Such a point exists because $S$ is rational.
The associated map is the stereographic
projection from $p$, which is birational to $\mathbb{P}^2$.
In this case, we have $P(-K_a)=2Q^2=4\le 6$.

Case (c): we take $P:=-K_a$, the class of hyperplane sections.
As parametrizing system, we choose the set of all sections with
hyperplanes containing the tangent plane through a fixed point $p$
defined over $\mathbb{K}$.
The associated map is the projection from the tangent plane,
which reduces the dimension by 3 and the degree by 4,
hence it is birational to $\mathbb{P}^2$. 
In this case, we have $P(-K_a)=P^2=5\le 6$.

Case (d): again we take $P:=-K_a$.
As parametrizing system, we choose the set of all sections with
hyperplanes containing the tangent plane through a fixed point $p$
and through another fixed point $q$ outside the tangent plane,
where both $p$ and $q$ are defined over $\mathbb{K}$.
The associated map can be decomposed into the projection 
from the tangent plane, which is birational onto a quadric in $\mathbb{P}^3$,
followed by the stereographic projection from the image of $q$.
In this case, we have $P(-K_a)=P^2=6$.
\end{proof}

\begin{theorem} \label{thm:i2}
Let $S$ be a nonsingular rational surface.
Let $P\in\mathrm{Pic}(S)$ be a free class
such that $P^2=0$ and $PK=-2$ and $\dim(P)=1$.
Assume that $S$ is $P$-minimal.
Then one of the following cases holds. \\
a) $S$ is isomorphic to the ruled surface $F_n$, $n\ge 0$. There exists
	an effective class $C$ such that $CP=1$, $C^2=-n$, and
	$K=(-n-2)P-2C$. The classes $C$ and $P$ generate $\mathrm{Pic}(S)$. \\
b) $S$ is isomorphic to the blowup of a nonsingular quadric 
	at a point of degree~2 (i.e. defined over a quadratic
	extension of $\mathbb{K}$). The Picard group is generated
	by $P$ and the exceptional class $E$, the class of
	plane sections is $P+E$, the canonical class is $-2P-E$,
	and $PE=2$. \\
c) $S$ is isomorphic to the blowup of $\mathbb{P}^2$ at a point
	of degree~4. The Picard group is generated by the
	exceptional class $E$ and the class of lines $L$, we have
	$P=2L-E$, and the canonical class is $-3L+E$.
\end{theorem}

\begin{proof}
Because $P$ is free, the associated map $S\to\mathbb{P}^1$ is regular,
and $P$ is the class of fibers. The genus of a generic fiber is
$\frac{P^2+PK}{2}+1=0$, hence the associated map gives $S$
the structur of a conic fibration. Let $d:=K^2$. Over the
algebraic closure $\overline{\mathbb{K}}$ is a blowup of a ruled surface $F_n$
for some $n$, hence $d\le 8$.

If $d=8$, then $S$ is minimal over $\overline{\mathbb{K}}$, and (a)
holds.

If $d=7$, then $S$ is not minimal by Theorem~4.1 in \cite{Iskovskih:70}.

If $d=6$, then (b) holds by Theorem~4.1 in \cite{Iskovskih:70}.

If $d=5$, then (c) holds by Theorem~4.1 in \cite{Iskovskih:70}.

If $d=4$, then $S$ cannot be both minimal and rational over $\mathbb{K}$,
by Theorem~2 in \cite{Iskovskih:72}.

If $d=3$, then $S$ cannot be both minimal and rational over $\mathbb{K}$,
by Corollary~2.6 in \cite{Iskovskih:70}.

If $d=2$ or $1$, then $S$ cannot be both minimal and rational over $\mathbb{K}$,
by Corollary~1.7 in \cite{Iskovskih:70}.

If $d\le 0$, then $S$ cannot be both minimal and rational over $\mathbb{K}$,
by Theorem~1.6 in \cite{Iskovskih:67}.
\end{proof}

\begin{theorem} \label{thm:j2}
Let $S$ be a nonsingular rational surface and let $D$ be a nef
divisor such that $D^2>0$ and $\mathrm{keel}(S)>0$. Then
\[ \pdeg(D) \le 4\ \mathrm{level}(D) + 2\ \mathrm{keel}(D) . \]
\end{theorem}

\begin{proof}
By lemmas~\ref{lem:c0}, \ref{lem:c1}, \ref{lem:c2}, and by replacing
$D$ by $2D$, we can reduce
to the case $\mathrm{level}(D)=a$ and $D_a=kP$ for $k=\mathrm{keel}(D)$
and $P$ as in Lemma~\ref{lem:c1}. Then we construct a parametrizing class $Q$
for each of the cases that arise in Theorem~\ref{thm:i2}. 
If $\phi_a:S\to S_a$ be the minimalization map,
then $\phi_a^\ast(Q)$ is a parametrizing class for $S$; and we
will prove the required upper bound for $\phi_a^\ast(Q)D$ in each case.

In case (a), we distinguish two subcases. If $S\cong F_0$, then
we take $Q:=C+P$. The image of the associated map is a ruled
quadric in $\mathbb{P}^3$. Because a stereographic projection
from a point defined over $\mathbb{K}$ is birational onto the plane,
the class $Q$ is parametrizing for $S_a$, and it follows that
$\phi_a^\ast(Q)$ is parametrizing for $S$. We compute
\[ \phi_a^\ast(Q)D = Q{\phi_a}_\ast(D) = Q(D_a-aK_a) \]
\[ = (C+P)(2aC+2aP+kP) = 4a+k \le 4a+2k. \]

If $S\cong F_n$ and $n>0$, then we take $Q:=C+nP$. 
Then $\dim(C+nP)=n+1$, and the image 
$S_Q\subset\mathbb{P}^{n+1}$ of the associated
map is a cone over a rational normal curve of degree~$n$. By repeating
projections from nonsingular points defined over $\mathbb{K}$,
we obtain a birational map to the plane. Hence $Q$ is parametrizing,
and therefore $\phi_a^\ast(Q)$ is parametrizing. We have
\[ \phi_a^\ast(Q)D = (C+nP)(2aC+(an+k+2a)P) = an+2a+k . \]
Moreover, $D$ is nef, hence 
\[ 0 \le \phi_a^\ast(D)(C) = 2a+k-an , \]
hence $\phi_a^\ast(Q)D\le 4a+2k$.

In case (b), we take $Q:=P+E$. The image of the associated map
is a quadric in $\mathbb{P}^3$, hence $Q$ is parametrizing, as
above. In this case, we get
\[ \phi_a^\ast(Q)D = (P+E)(kP+2aP+aE) = 4a+2k . \]

In case (c), we take $Q=L$, which is of course parametrizing.
In this case, we get
\[ \phi_a^\ast(Q)D = L(kP-aK) = L(3aL+2kL-(a+k)E) = 3a+2k \le 4a+2k . \]
\end{proof}

Comprizing Theorems~\ref{thm:lb}, \ref{thm:j1}, \ref{thm:j2}, and 
Remark~\ref{rem}, we can finally state:

\begin{theorem} \label{thm:fin}
For any rational surface $S$, the following bounds hold:
\[ 3\ \mathrm{level}(S)+\mathrm{keel}(S) \le \pdeg(S)
\le 6\ \mathrm{level}(S)+2\ \mathrm{keel}(S) . \]
\end{theorem}

\bibliography{all}
\bibliographystyle{plain}
 
\end{document}